\documentclass[a4paper]{amsart}

\usepackage{amsmath, amsthm, amssymb}
\usepackage{eucal, enumerate, cite, graphicx}
\usepackage[all]{xy}

\theoremstyle{plain}
\newtheorem{theorem}{Theorem}[section]
\newtheorem{proposition}[theorem]{Proposition}
\newtheorem{corollary}[theorem]{Corollary}

\theoremstyle{definition}
\newtheorem{definition}[theorem]{Definition}
\newtheorem{conjecture}[theorem]{Conjecture}
\newtheorem{example}[theorem]{Example}

\theoremstyle{remark}
\newtheorem{remark}[theorem]{Remark}

\begin{document}
\bibliographystyle{alpha}

\title[]{About the hyperbolicity of complete intersections}

\author{Simone Diverio}
\address{Simone Diverio \\ CNRS and Institut de Math\'ematiques de Jussieu - Paris Rive Gauche}
\email{diverio@math.jussieu.fr} 
\thanks{The author is partially supported by the ANR project \lq\lq POSITIVE\rq\rq{}, ANR-2010-BLAN-0119-01.}

%\date{\today}

\begin{abstract}
This note is an extended version of a thirty minutes talk given at the \lq\lq XIX Congresso dell'Unione Matematica Italiana\rq\rq{}, held in Bologna from September 12th to September 17th, 2011. This was essentially a survey talk about connections between Kobayashi hyperbolicity properties and positivity properties of the canonical bundle of projective algebraic varieties. 
\end{abstract}

\maketitle

\section{Notations and preliminary material}

Let $X$ be a compact complex space. 

\begin{definition}
An \emph{entire curve} in $X$ is a non constant holomorphic map $f\colon\mathbb C\to X$.
\end{definition}

We shall give a definition of Kobayashi hyperbolicity which is often referred to as Brody hyperbolicity. Nevertheless, by Brody's lemma, the two concepts coincide in the compact case: in the sequel, we shall freely speak of \lq\lq hyperbolicity\rq\rq{} referring to these equivalent concepts. For more details on Brody's lemma and Brody hyperbolicity, see for instance \cite{sd_Kob98}.

\begin{definition}
A compact complex space $X$ is \emph{hyperbolic} if and only if there are no entire curves in $X$.
\end{definition}

\begin{remark}
Fix any hermitian metric $\omega$ on $X$. By Brody's lemma, if there exists on $X$ an entire curve, then there exists on $X$ also a so-called Brody curve, \textsl{i.e.} an entire curve $f$ with bounded derivative $||f'||_\omega\le C$. 
\end{remark}

Let us give some basic examples of hyperbolic compact manifolds.

\begin{example}
Let $C$ be a compact Riemann surface. Of course if $C$ is the Riemann sphere or an elliptic curve, then $C$ is not hyperbolic. On the other hand, by the Uniformization Theorem, if the geometric genus of $C$ is greater than or equal to two, then the universal cover of $C$ is the complex unit disc. Thus, since every holomorphic map $f\colon\mathbb C\to C$ lifts to a holomorphic map to the universal cover, by Liouville's theorem it must be constant. Therefore, $C$ is hyperbolic if and only if $g(C)\ge 2$.
\end{example}

\begin{example}[Kobayashi]
Let $X$ be a compact complex manifold with ample cotangent bundle $T^*_X$, \textsl{i.e.} such that there exists a positive integer $k_0$ such that all evaluation maps 
$$
H^0(X,S^kT^*_X)\to J^1S^kT^*_{X,x},\quad H^0(X,S^kT^*_X)\to S^kT^*_{X,x}\oplus S^kT^*_{X,y},
$$
for $x,y\in X$, $x\ne y$, are surjective for all $k\ge k_0$. Here, $S^kT^*_X$ is the $k$-th symmetric power of the cotangent bundle and $J^1S^kT^*_X$ is the bundle of $1$-jets of sections of $S^kT^*_X$. 

We claim that $X$ is hyperbolic. To see this, fix an integer $k\ge k_0$ as above and a basis $\sigma_0,\dots,\sigma_N$ of $H^0(X,S^kT^*_X)$. Then, the function
$$
\begin{aligned}
\eta\colon & T_X\to\mathbb R \\
 & v \mapsto \biggl(\sum_{j=0}^N |\sigma_j(x)\cdot v^{\otimes k}|^2\biggr)^{1/2k},\quad v\in T_{X,x},
\end{aligned}
$$
is easily seen to be plurisubharmonic, and strictly plurisubharmonic outside the zero section, thanks to the ampleness assumption. Now, suppose that $X$ is not hyperbolic. Then, there exists on $X$ a non constant Brody curve $f\colon\mathbb C\to X$ and therefore $\eta\circ f\colon\mathbb C\to\mathbb R$ is a bounded subharmonic function defined on the entire complex plane which is strictly subharmonic outside $\{f'=0\}$. Hence, $f$ is constant and we get a contradiction.

Examples of compact complex manifolds with ample cotangent bundle are given for instance by the intersection of at least $n/2$ sufficiently ample general hypersurfaces in an abelian variety of dimension $n$, and by general linear section of small dimension in a product of sufficiently many smooth projective varieties with big cotangent bundle (see \cite{sd_Deb05} for more details). Other examples are given by compact hermitian manifolds with negative holomorphic bisectional curvature.
\end{example}

It is conjectured in \cite{sd_Deb05} that the intersection, in $\mathbb P^N$, of at least $N/2$ generic hypersurfaces of sufficiently high degrees has ample cotangent bundle (and hence is hyperbolic). We shall come back to that later.

\begin{example}
More generally, a compact hermitian manifold $(X,\omega)$ such that its holomorphic sectional curvature is negative is hyperbolic. Any compact quotient $\Delta^n/\Gamma$, $n\ge 2$, of a polydisc by a group $\Gamma\subset\operatorname{Aut}(\Delta)^n$ acting freely and properly discontinuously on $\Delta^n$ gives an example of compact hermitian manifold with negative holomorphic sectional curvature but with no hermitian metrics of negative holomorphic bisectional curvature. 
\end{example}

In the case of complex projective curves, the hyperbolicity can be characterized in a completely algebraic way: a complex projective curve is hyperbolic if and only if its geometric genus is greater than or equal to two, and this happens if and only if its cotangent bundle, which in complex dimension one coincides with its canonical bundle, is ample. By the Riemann-Roch formula, this is equivalent to being of general type. Let us give the precise definition.

\begin{definition}
Let $X$ be a compact complex manifold and $L\to X$ a holomorphic line bundle. The \emph{Kodaira-Iitaka dimension $\kappa(L)$ of $L$} is defined as follows: if $H^0(X,L^{\otimes m})=\{0\}$ for each $m\ge 1$ then put $\kappa(L)=-\infty$, otherwise $0\le \kappa(L)\le\dim X$ is the unique integer such that 
$$
C^{-1}\,m^{\kappa(L)}\le \dim H^0(X,L^{\otimes m}) \le C\,m^{\kappa(L)},
$$ 
for some positive constant $C$ and for all sufficiently large and divisible integers $m$.
If $\kappa(L)=\dim X$, $L$ is said to be \emph{big}. We denote by $\kappa(X)$ the Kodaira-Iitaka dimension of the canonical bundle $K_X=\det T^*_X$ of $X$ and call it the \emph{Kodaira dimension of $X$}. The manifold $X$ is said to be of \emph{general type} if it has maximal Kodaira dimension.
\end{definition}

\begin{remark}
By Hirzebruch-Riemann-Roch and Kodaira's vanishing, if $L\to X$ is ample then $L$ is big. Thus, manifolds with ample canonical bundle are of general type. Beside projective curve of genus $\ge 2$, a typical example of manifold of general type is given by smooth projective hypersurfaces $X\subset\mathbb P^{n+1}$ such that $\deg X\ge n+3$ (they have in fact ample canonical bundle, as it can be straightforwardly checked by adjunction).
\end{remark}

\section{Guiding conjectures}

Let us begin with the following proposition.

\begin{proposition}[Demailly \cite{sd_Dem97}]
Let $X$ be a compact complex manifold and $\omega$ a hermitian form on $X$. Consider the following statements:
\begin{enumerate}[(i)]
\item $X$ is Kobayashi hyperbolic.
\item There exists $\varepsilon_0>0$ such that for every curve $C\subset X$ 
$$
-\chi(\widehat C)=2g(\widehat C)-2\ge\varepsilon_0\,\deg_\omega C,
$$
where $\widehat C\to C$ is the normalization of $C$ and $\deg_\omega C:=\int_C\omega$.
\item Every holomorphic map $Z\to X$, where $Z$ is a complex torus, is constant.
\end{enumerate}
Then, $(i)\implies (ii)\implies (iii)$.
\end{proposition}

Property $(ii)$ above is often referred to as \emph{algebraic hyperbolicity}. Now, if $X$ is supposed to be projective algebraic, the three properties above are conjectured to be equivalent.

\begin{conjecture}[Lang \cite{sd_Lan86}, Demailly \cite{sd_Dem97}]\label{langdem}
In the proposition above, if $X$ is moreover projective algebraic, then $(iii)\implies (i)$.
\end{conjecture}

\begin{remark}
If $X$ is supposed to be merely K\"ahler, then the conjecture is false. One can show, for instance, that a $K3$ surface without any curve does not admit any non constant map from a complex torus. Nevertheless, a $K3$ surface is never hyperbolic. For more details, see \cite{sd_Can00}.
\end{remark}

Thus, from now on, $X$ will be always projective algebraic. In this case, it is very tempting to believe that the analytic property of being hyperbolic should translate to completely algebraic properties of the manifold (generalizing then the case of complex dimension one). This belief is formalized in the following.

\begin{conjecture}[Lang \cite{sd_Lan86}]\label{lang}
Let $X$ be a complex projective manifold. Then, $X$ is Kobayashi hyperbolic if and only if $X$ as well as all its subvarieties are of general type.
\end{conjecture}

The sufficiency of the condition in the conjecture above is the object of this other one.

\begin{conjecture}[Green-Griffiths]
Let $X$ be a projective manifold of general type. Then, there should exists a proper algebraic subvariety $Y\subsetneq X$ such that all entire curves in $X$ are actually contained in $Y$.
\end{conjecture}

The necessity of the condition in Conjecture \ref{lang}, since a projective manifold of general type without rational curves has ample canonical bundle, would imply the following.

\begin{conjecture}[Kobayashi]
Let $X$ be a projective hyperbolic manifold. Then, the canonical bundle $K_X$ should be ample.
\end{conjecture}

Observe that all these conjectures are obvious in dimension one, but already non trivial (and, apart from the last one, not known) in dimension two. 

To end this section, let us mention that there is also an arithmetic counterpart of the story, which is moreover already highly non trivial in dimension one. In number theory, the Mordell conjecture states that a curve of genus greater than $1$ defined over a number field $K$ has only finitely many rational points over any finite field extension $L$ of $K$: it was proved by G. Faltings in 1983. Remind that until Faltings' theorem, there was not known a single example of a curve which was proved to have only a finite number of rational points in every finitely generated field over $\mathbb Q$. \lq\lq Accidentally\rq\rq{}, the curves satisfying Faltings' theorem are exactly the hyperbolic ones. According to Lang, this should be the right framework in order to have a higher dimensional analogous of Faltings' theorem.

\begin{conjecture}[Lang \cite{sd_Lan74}]
Let $X$ be a projective manifold defined over a number field $K$. Then, $X$ is hyperbolic if and only if it contains only finitely many rational points over any finite field extension $L$ of $K$.
\end{conjecture}

At present, almost nothing is known about this very fascinating theme.

\section{Overview of old and recent results}

We now survey some results relative to the conjectures mentioned above.

\subsection{The Kobayashi conjecture}

As mentioned in the previous sections, for the case of curves this conjecture is a trivial consequence of the Uniformization Theorem and Liouville's theorem. 

For surfaces, the Enriques-Kodaira birational classification implies that, in order to prove the conjecture, it suffices to show that K3 surfaces are not hyperbolic. Since Kummer surfaces (which are obviously non hyperbolic) are dense in the moduli space of K3's and since hyperbolicity is an open property (with respect to analytic topology), the conjecture follows. In fact, much more is true: by \cite{sd_MM83} every projective K3 surface does contain at least one rational curve. 

Unfortunately, we do not dispose of an analogous result in higher dimension. It is likely that every compact projective manifold with vanishing first real Chern class is not hyperbolic and even more: it should admit a non trivial holomorphic map from a complex torus (or, more optimistically, a rational curve).

Now, we concentrate on dimension three. To begin with, observe that several powerful machineries from birational geometry ---such as the characterization of uniruledness in terms of negativity of the Kodaira dimension (which holds true up to dimension three and is conjectural in general), the Iitaka fibration, the abundance conjecture (which is actually a theorem up to dimension three)--- permit to reduce this conjecture to the following statement: a projective threefold $X$ of Kodaira dimension $\kappa(X)=0$ cannot be hyperbolic. This is the analogous of the reduction to K3 surfaces in dimension two. By abundance, the Beauville-Bogomolov decomposition theorem and elementary properties of hyperbolic manifolds, in order to prove the Kobayashi conjecture in dimension three \emph{it suffices to show that a Calabi-Yau threefold is not hyperbolic}. Here, by a Calabi-Yau threefold we mean a simply connected compact projective threefold with trivial canonical class $K_X\simeq\mathcal O_X$ and $H^i(X,\mathcal O_X)=0$, $i=1,2$.

All known results at present prove in fact a much stronger statement, namely existence of rational curves on Calabi-Yau threefolds. Let us cite a couple of results in this direction:

\begin{itemize}
\item The article \cite{sd_HB-W92} is the culmination of a series of papers by Wilson in which he studies in a systematic way the geometry of Calabi-Yau threefolds; among many other things, it is shown there that if the Picard number $\rho(X)>13$, then there always exists a rational curve on $X$.
\item Following somehow the same circle of ideas, it was proven in \cite{sd_Pet91} (see also \cite{sd_Ogu93}) that a Calabi-Yau threefold $X$ has a rational curve provided there exists on $X$ a non-zero effective non-ample line bundle on $X$.
\end{itemize}

By the Cone Theorem, if there exists on a Calabi-Yau manifold $X$ a non-zero effective non-nef line bundle, then there exists on $X$ a rational curve (generating an extremal ray). Therefore, we can always suppose that such an effective line bundle is nef. Remark, on the other hand, that in Peternell's result, the effectivity hypothesis is crucial (regarding it in a more modern way) in order to make the machinery of the logMMP work. In this spirit, Oguiso asked in \cite{sd_Ogu93} the following question: \emph{is it true that if a Calabi-Yau threefold $X$ possesses a non-zero nef non-ample line bundle, then there exists a rational curve on $X$?}

Here is a positive answer, under a mild condition on the Picard number of $X$.

\begin{theorem}[Diverio-Ferretti \cite{sd_D-F}]
Let $X$ be a Calabi-Yau threefold and $L\to X$ a non-zero nef non-ample line bundle. Then, $X$ has a rational curve provided $\rho(X)>4$.
\end{theorem}

The proof of this theorem is of arithmetic nature and relies upon a careful study of the diophantine properties of the cubic intersection form defining the nef boundary in the N\'eron-Severi space of $X$. In particular, it is exploited the link between the boundedness and convexity properties of the nef cone and the arithmetic of its boundary. 

The results stated above thus permit, as in \cite{sd_Pet91}, to exclude a certain number of cases to be checked in order to prove such a non-hyperbolicity statement. Finally, as far as we know, there is no known example of a Calabi-Yau threefold without non-constant holomorphic images of complex tori.

\subsection{The Green-Griffiths-Lang conjecture}

This conjecture in full generality is still unknown even in dimension two, so let us start by discussing this case. 

To begin with, observe that if the Zariski closure of all entire curves in a surface $S$ is positive dimensional (otherwise the surface is hyperbolic) but a proper subvariety, then by uniformization it must be a finite union of rational and elliptic curves.
If the second Segre number $c_1^2-c_2$ of a minimal surface of general type $S$ is positive, then it is know, since the work of Bogomolov, that there is abundance of symmetric differentials on $S$. Thanks to this, Bogomolov himself was able to show that there is only a finite number of rational and elliptic curves on such a surface $S$, see \cite{sd_Des79}. 

Several years later, McQuillan \cite{sd_McQ98}, proving deep results on holomorphic (multi)\-foliations on surfaces, was able to extend the work of Bogomolov to the transcendental setting, thus showing that the Green-Griffiths conjecture is true for surfaces of general type with positive second Segre number. More precisely, using Miyaoka's semi-positivity result for cotangent bundles of nonuniruled projective varieties and a dynamic diophantine approximation, he derived strong Nevanlinna Second Main Theorems for entire curves tangent to the leaves of a holomorphic foliation. In particular, he obtains that every parabolic leaf of an algebraic (multi)foliation on a surface $S$ of general type is algebraically degenerate. The assumption $c_1^2 > c_2$ guarantees the existence of an algebraic multi-foliation such that every entire curve is contained in one of its leaves. 

This last sentence can be rephrased by saying that under the assumption on the second Segre class, every entire curve must satisfy an algebraic differential equation of order one.
This can be put in perspective as follows (see \cite{sd_Dem97} for all the details). Let $X$ be a complex manifold, $p_k\colon J_kX\to X$ the holomorphic fiber bundle of $k$-jets of germs of holomorphic curves $\gamma\colon(\mathbb C,0)\to X$ and $J_kX^{\text{reg}}$ its subset of regular ones, \textsl{i.e.} such that $\gamma'(0)\ne 0$. There is a natural action of the group $\mathbb G_k$ of $k$-jets of biholomorphisms of $(\mathbb C,0)$ on $J_kX$, and the quotient $J_kX^{\text{reg}}/\mathbb G_k$ admits a nice geometric relative compactification $J_kX^{\text{reg}}/\mathbb G_k\hookrightarrow X_k$. Here, $\pi_k\colon X_k\to X$ is a tower of projective bundles over $X$. In particular, it is naturally endowed with a tautological line bundle $\mathcal O_{X_k}(-1)$.

It turns out that the direct image $(\pi_k)_*\mathcal O_{X_k}(m)$, $m\ge 1$, is the sheaf of holomorphic sections of a holomorphic vector bundle $E_{k,m}T^*_X\to X$ on $X$. For instance, when $k=1$, we have that $E_{1,m}T^*_X\simeq S^mT^*_X$ is just the $m$-th symmetric power of the cotangent bundle. 

\begin{definition}
The holomorphic sections of the vector bundle $E_{k,m}T^*_X$ are called \emph{invariant jet differentials of order $k$ and (weighted) degree $m$}. 
\end{definition}

Invariant jet differentials of order $k$ and (weighted) degree $m$ act on ($k$-jets of) holomorphic curves traced in $X$ as polynomial differential operators. More precisely, let $P\in\mathcal O_X(E_{k,m}T^*_X)(U)$ be an invariant jet differentials over an open set $U\subset X$. Then,  $P$ is nothing but a holomorphic function $P\colon p_k^{-1}(U)\subset J_kX\to\mathbb C$, such that $P(\gamma\circ\varphi)=(\varphi')^m P(\gamma)$, for $\gamma\colon(\mathbb C,0)\to U$ a $k$-jet of germs of holomorphic curve and $\varphi\colon(\mathbb C,0)\to (\mathbb C,0)$ a $k$-jet of biholomorphism of the origin. This immediately explains why $E_{1,m}T^*_X\simeq S^mT^*_X$ and why, in general, invariant jet differentials can be considered as algebraic differential operators acting on holomorphic curves.

\begin{theorem}[Green-Griffiths \cite{sd_GG80}, Siu-Yeung \cite{sd_SY96}, Demailly \cite{sd_Dem97}] \label{}
Let $X$ be projective algebraic and $A\to X$ an ample line bundle. Then, for all invariant jet differentials $P\in H^0(X,E_{k,m}T^*_X\otimes A^{-1})\simeq H^0(X_k,\mathcal O_{X_k}(m)\otimes\pi_k^* A^{-1})$ with values in the anti-ample divisor $A^{-1}$ and all entire curves $f\colon\mathbb C\to X$, we have $P(f)\equiv 0$.
\end{theorem}

Therefore, abundance of invariant jet differentials is translated by this theorem in abundance of constraints for entire curves. From this point of view, the first step to prove hyperbolicity type (or, more generally, algebraic degeneracy type) statements is to produce some invariant jet differential. 

For instance, it is conjectured (and proved for surfaces) by Green and Griffiths in \cite{sd_GG80} that if $X$ is a projective algebraic manifold of general type then there should exist a $k\gg 1$ such that $\mathcal O_{X_k}(1)$ is big (or rather, another less refined version of this line bundle is big, when $X_k$ is obtained by modding out just by the group of homotheties instead of the full group of reparametrization of the origin). 
This conjecture has been proved in full generality only very recently by J.-P. Demailly in \cite{sd_Dem11}, by means of  an astonishing combination of his holomorphic Morse inequalities and a probabilistic interpretation of higher order jets.

\subsubsection{Projective hypersurfaces and complete intersections} 

Since in full generality the Green-Griffiths-Lang conjecture seems out of reach for the moment, several efforts have been made in some particular cases. Let us concentrate on projective hypersurfaces and, later on, more generally, on complete intersections.

Let $X\subset\mathbb P^{n+1}$, $n\ge 2$, be a projective hypersurface of degree $\deg X=d$. If $X$ is smooth, then adjunction formula together with the Euler exact sequence show that $K_X\simeq\mathcal O_X(d-n-2)$. Thus, $K_X$ is ample if (and only if) $d\ge n+3$; in particular, if $\deg X\ge n+3$, then $X$ is of general type. Therefore, if $X$ is a smooth projective hypersurface of degree $\deg X\ge n+3$, the Green-Griffiths conjecture predicts the existence of a proper subvariety $Y\subsetneq X$ containing all entire curve in $X$.

If, moreover, $\deg X\ge 2n+1$ and $X$ is (very) generic, then Voisin proved in \cite{sd_Voi96,sd_Voi96err} (among other things) that all its subvarieties are of general type (this was first proven in \cite{sd_Ein88}) and, as a consequence of her proof, that $X$ is algebraically hyperbolic. Here, by very generic we mean a hypersurface whose modulus is outside a countable union of proper subvarieties of the parameter space. Then, in view of Conjectures \ref{langdem} and \ref{lang}, in this case $X$ should be hyperbolic: this was in fact the content of the following conjecture made by Kobayashi in 1970.

\begin{conjecture}[Kobayashi]
Let $X\subset\mathbb P^{n+1}$, $n\ge 2$, be a smooth (very) generic projective hypersurface of degree $\deg X\ge 2n+1$. Then, $X$ is hyperbolic.
\end{conjecture}

Observe, first of all, that the genericity assumption is necessary since, for example, Fermat's hypersurfaces always contain lines (which shows also that ample canonical bundle does not imply hyperbolic). Second, the conjecture implies in particular that every generic projective complete intersection of high (multi)degree should be hyperbolic.

For surfaces in projective $3$-space, the conjecture is known without optimality of the degree $d$, which should be starting from five: it has been proven independently by McQuillan \cite{sd_McQ99} for $d\ge36$ and Demailly-El Goul \cite{sd_DEG00} for $d\ge 21$. The lower bound on the degree has been further refined by P\u{a}un \cite{sd_Pau08} to $d\ge 18$.

For threefolds in projective $4$-space, Rousseau \cite{sd_Rou07} was able to show that if $\deg X\ge 593$ and $X$ is generic, then there does not exist any Zariski dense entire curve in $X$ (but observe that the locus covered by the image of entire curve is not excluded to be Zariski dense).

Next, let us give a brief account about more recent and somehow general results.

\begin{theorem}[Diverio-Merker-Rousseau \cite{sd_DMR10}, Diverio-Trapani \cite{sd_DT11}]\label{dmr}
Let $X\subset\mathbb P^{n+1}$ be a generic projective hypersurface of degree $\deg X\ge 2^{n^5}$. Then, there exists a proper algebraic subvariety $Y\subsetneq X$ which contains the image of all entire curves in $X$. If $\dim X\ge 3$ or $\dim X=2$ and $X$ is very generic, then the codimension of $Y$ in $X$ is at least two.
\end{theorem}

The information about the codimension of $Y$ is the contribution of \cite{sd_DT11}: in order to obtain it, it is crucial that every effective divisor is ample (this explains the very genericity assumption for surfaces, in view of the Noether-Lefschetz theorem).

Before discussing the key ingredients in the proof of the theorem above, let us make some comments and state a couple of corollaries. As far as we know, this is the first result in all dimensions about algebraic degeneracy of entire curves in manifolds of general type. On the other hand, unfortunately, projective hypersurfaces form a too small class to see this result really as an evidence for the Green-Griffiths conjecture. Nevertheless, this theorem gives also an effective estimate on the degree of the hypersurfaces only in terms of their dimension (even if this degree is very far from being optimal...). It would be desirable to drop the genericity assumption (which is intrinsic in the methods used in the proof) as far as only algebraic degeneracy (and not the full hyperbolicity) is concerned.

\begin{corollary}[Diverio-Trapani \cite{sd_DT11}]
A very generic threefold in $\mathbb P^4$ of degree at least $593$ is Kobayashi hyperbolic.
\end{corollary}

This follows from the fact that the degeneracy locus is of dimension one in this case and thus it should consist of rational or elliptic curves, which are absent in these hypersurfaces.

\begin{corollary}[Brotbek \cite{sd_Bro11}]
Let $X$ be a generic projective complete intersection of (multi)degree greater than or equal to $(2^{n^5},\dots,2^{n^5})$. Then, $X$ is hyperbolic provided $\dim X\le 2\operatorname{codim} X$.
\end{corollary}

This corollary is obtained by one further genericity argument combined with the action of the group of projective automorphisms. Observe that this last result can be seen as a partial confirmation of Debarre's conjecture stated in the introduction (for more developments on Debarre's conjecture, especially for generic complete intersection surfaces of high degree in $\mathbb P^4$ which are shown to have ample cotangent bundle, see \cite{sd_Bro11}).

\begin{proof}[Heuristic idea of the proof of Theorem \ref{dmr}]
Thanks to the strategy proposed in the work \cite{sd_Siu04} of Siu and inherited somehow by Voisin's work \cite{sd_Voi96,sd_Voi96err}, it was quite clear that in order to prove such an algebraic degeneracy statement one had, above all, to overcome two main difficulties:

\begin{itemize}
\item[(i)] Produce some invariant jet differential vanishing on an ample divisor on every smooth projective hypersurface of large degree.
\item[(ii)] Produce low pole order meromorphic vector fields on the universal hypersurface, in order to produce by differentiation new jet differentials, enough to control the geometry of entire curves.
\end{itemize}

Step (i) has been done in \cite{sd_Div09}, as an application of algebraic holomorphic Morse inequalities: the method is intrinsically effective and permits to estimate explicitly, in terms of the dimension only, how large the degree should be. Step (ii) has been achieved in \cite{sd_Mer09}. 

Once we have at our disposal a \lq\lq first\rq\rq{} jet differentials, thanks to a semi-continuity argument we can extend it holomorphically to all projective hypersurfaces parame\-trized by a Zariski open set of their moduli space (this semi-continuity argument makes the result \lq\lq only\rq\rq{} generic). Next, the extended jet differential can be seen as a jet differential on an open set of the universal hypersurface. Therefore, one can now use low pole order meromorphic tangent vector fields to derive this algebraic operator in order to obtain new ones, algebraically independent from the first one. These new operators, once restricted to every projective hypersurface parametrized by this Zariski open set, give enough informations to conclude that their base locus projects onto a proper subvariety. The low pole order property enables us to be sure that after differentiation the jet differential remains holomorphic and with values in an anti-ample divisor.
\end{proof}

Very recently, Siu, in \cite{sd_Siu12}, was able to push forward his own strategy in order to obtain Kobayashi hyperbolicity of generic projective hypersurfaces of high degree. The details of his beautiful and delicate proof are underway of validation by the experts.

\end{document}